\documentclass{article}

\usepackage{amssymb}

\usepackage{epsfig}

\newcommand{\altheta}{\alpha(\theta)}

\newcommand{\detd}{\det \,  d (\pi \circ  \phi _{t}) _{\theta}}

\newcommand{\detdt}{\det \,  d (\pi \circ  \phi _{t}) _{(\theta,t)}}

\newcommand{\htop}{h_{top}}

\newcommand{\entreceroyt}{\int_{0}^{T} dt}

\newcommand{\enSM}{\int_{SM}}

\newcommand{\nst}{n_{T}(x,y)}

\newcommand{\reals}{\mathbb{R}}

\newtheorem{lemma}{\bf{Lemma}}

\newtheorem{teo}{\bf{Theorem}}

\title{Topological entropy of a magnetic flow and the growth of the number of trajectories}

\author{C\'esar J. Niche \thanks{This work was partially supported by the NSF. The author acknowledges support from IMERL, Fac. de Ingenier\'{\i}a - UDELAR, Montevideo, 
Uruguay.}\\ Mathematics Department \\ University of California at Santa Cruz \\ Santa Cruz CA 95064 \\ USA \\ E-mail: cniche@math.ucsc.edu}

\begin{document}

\noindent

\maketitle

\begin{abstract}
We prove formulae relating the topological entropy of a magnetic flow to the growth rate of the average number of trajectories connecting two points.
\end{abstract}

\section{Introduction}

Topological entropy measures how complicated the global orbit structure of a dynamical system is, in terms of the exponential growth rate of the number of orbit segments that can be distinguished with arbitrarily fine, but finite, precision. It is important to relate this conjugacy invariant to dynamical or geometrical characteristic elements of the specific type of dynamical system under study, in order to have more tools to understand the structure of such systems. \\
For the geodesic flow of a $C^{\infty}$ metric on a manifold $M$, Ma\~n\'e \cite{manhe} proved formulae relating the topological entropy of the flow to the growth rate of: 

\begin{itemize}
\item[(A)] the average on the unit sphere bundle of $M$ of the determinant of the differential of the flow on the vertical distribution $V(\theta) = \ker (d\pi)_{\theta}$, for $\pi: SM \to M$ the usual projection

\item[(B)]  the average of the number of geodesics connecting any two points on the manifold $M$.
\end{itemize}

Ma\~n\'e's proof of this Riemannian result is symplectic, as the geodesic flow is Hamiltonian for the kinetic energy on $TM$ endowed with a natural symplectic structure. One of the key points in the proof is the fact that the Lagrangian vertical distribution $\bar{V}(\theta) = \ker (d\bar{\pi})_{\theta}$, where $\bar{\pi}: TM \to M$ is the standard projection, ``is always twisted by the flow in the same direction'' (see Paternain \cite{pat}). \\

This twist property of the geodesic flow for the vertical distribution, is an instance of the optical property of a Hamiltonian flow for a Lagrangian distribution, introduced by Bialy and Polterovich \cite{bp}. For optical Hamiltonian flows, an appropiate generalization of (A) was proved by Niche \cite{niche3}. \\

An important family of optical Hamiltonian flows is that of magnetic flows (or twisted geodesic flows). Let $(M, \langle \cdot, \cdot \rangle)$ be a smooth closed Riemannian manifold and $\Omega$ a closed $2$-form on $M$. The tangent bundle $TM$ is  a symplectic manifold with symplectic structure $\omega = \omega_{L} + \pi ^{\ast} \Omega$, where $\omega_{L}$ is the pullback of the standard symplectic form of $T^{\ast}M$ by the metric. The magnetic flow $\phi _{t}$ is the Hamiltonian flow of $H:TM \to \reals$, $H(x,v) = \frac{1}{2} \langle v,v \rangle _{x}$ with respect to the form $\omega$. \\

The result in this note is the following version of (B) for magnetic flows.

\begin{teo}
\label{main-teo}
Let $\phi _{t}$ be a $C ^{\infty}$ magnetic flow and let $\nst$ be the number of trajectories of the flow connecting  $x$ and $y$ with length less than or equal to $T$. Then

\begin{displaymath}
\htop \geq \limsup _{T \to + \infty} \frac{1}{T} \log \int _{M \times M} \nst \, dx \, dy
\end{displaymath}

where $dx \, dy$ is the product Riemannian measure in $M \times M$.

If the flow admits a codimension one, continuous, invariant distribution of hyperplanes in $SM$ transversal to the magnetic vector field, then
\begin{displaymath}
\htop = \lim _{T \to + \infty} \frac{1}{T} \log \int _{M \times M} \nst \, dx \, dy.
\end{displaymath}
\end{teo}

{\bf Remarks}

\begin{itemize}
\item The sum of the strong stable and strong unstable bundles of an Anosov flow provides a H\"older continuous, codimension one, invariant distribution of hyperplanes, transversal to the vector field, so our second formula applies to magnetic Anosov flows.

\item Another important class of magnetic flows which admit a codimension one, transversal, invariant distribution is that of contact type magnetic flows. In dimension two and for $M \neq T^{2}$, the magnetic flow has contact type for high energy levels and, provided that the magnetic field $\Omega$ is symplectic, for low energy levels. In dimensions greater than two  the magnetic flow is contact type for high energy values if and only if the magnetic field $\Omega$ is exact. The distributions associated to these contact forms are $C^{\infty}$.

\item Applications of Ma\~ne's result can be found in Bolotin and Negrini \cite{bol-neg-1}, \cite{bol-neg-2}.

\item Related results about the entropy of magnetic flows can be found in Burns and Paternain \cite{burnspat}, Grognet \cite{grognet-1} and Paternain and Paternain \cite{pat-pat-2}. 

\end{itemize}

The proof of Theorem \ref{main-teo} is based on Theorem 5.1 from Burns and Paternain \cite{burnspat} and the following

\begin{lemma}
\label{main-lemma}
Let $\phi_{t}$ be a $C^{\infty}$ magnetic flow, $\pi: SM \to M$ the standard projection, $V(\theta) = \ker (d\pi)_{\theta}$ and $\altheta = V(\theta)  \oplus \langle X(\theta) \rangle$, where $X(\theta)$ is the magnetic vector field. Then

\begin{displaymath}
\int _{M \times M} n_{T}(x,y) \, dx \, dy =  \entreceroyt \enSM \vert \detd \vert_{\altheta} \vert \, d\theta.
\end{displaymath}
\end{lemma}

\section{Setting and proofs}

As we stated in the introduction, the magnetic flow is the Hamiltonian flow of the kinetic energy on the symplectic manifold $(TM, \omega = \omega _{L} + \pi ^{\ast} \Omega)$, for $\omega_{L}$ the pullback of the canonical symplectic form of $T^{\ast}M$ by the metric and $\Omega$ a closed $2$-form. The Lorentz force is the linear map $Y: TM \to TM$ such that

\begin{displaymath}
\Omega _{x} (u,v) = \langle Y(u), v \rangle, \quad u, v \in T_{x}M.
\end{displaymath}

Then the magnetic vector field $X$ can be expressed in $TTM$ as

\begin{displaymath}
X (\theta) = (v, Y(v))
\end{displaymath}

for $\theta = (x,v) \in TM$. Here we use the standard splitting $T_{\theta}TM = H(\theta) \oplus \bar{V}(\theta)$, where $H(\theta)$ is the lift of $T_{x}M$ by the metric and $\bar{V}(\theta)$ is the kernel of $(d\bar{\pi}) _{\theta}$, for $\bar{\pi}:TM \to M$. A curve $(\gamma, \dot{\gamma})$ in $TM$ is an integral curve of $X$ if and only if

\begin{displaymath}
\frac{D}{dt} \dot{\gamma} = Y(\dot{\gamma})
\end{displaymath}

for the covariant derivative $D$. This equation is just Newton's law for a charged particle of unit mass and charge. It is easy to check that the magnetic flow leaves $SM$ invariant. \\

{\em Proof of Lemma \ref{main-lemma}:} Let $f:S_{x}M \times [0,T] \to M$ be $f(\theta, t) = (\pi \circ \phi _{t}) (\theta)$. It is easy to see that $f^{-1}(y) = n_{T}(x,y)$, where $n_{T}(x,y)$ is the number of trajectories of the flow connecting $x$ and $y$ with length less or equal than $T$. The Area Formula (see page 243, Federer \cite{fed}) implies that

\begin{equation}
\label{eqn:first}
\int _{M} n_{T}(x,y) \, dy = \int _{S _{x}M \times [0,T]} \vert \detdt \vert_{T_{(\theta,t)}(S_{x} M\times [0,T])} \vert  \, d(\theta,t).
\end{equation}

Identifying the spaces as follows


\begin{displaymath}
T_{(\theta,t)}(S_{x} M\times [0,T]) = T_{\theta}(S_{x} M) \oplus T_{t}([0,T]) = V(\theta) \oplus \langle X(\theta) \rangle = \altheta
\end{displaymath}

we see that the maps 

\begin{displaymath}
d (\pi \circ  \phi _{t}) _{(\theta,t)}: T_{(\theta,t)}(S_{x} M\times [0,T]) \to T _x M,  \quad d (\pi \circ  \phi _{t}) _{(\theta,t)}: \altheta \to T _x M
\end{displaymath}

are equal. Let us endow $S_x M \times [0,T]$ with the product metric. Note that in this metric, $V(\theta)$ and $X(\theta)$ are orthogonal, and hence $d(\theta,t) = d\theta \wedge dt$, while $V(\theta)$ and $X(\theta)$ are not orthogonal with respect to the metric induced on $TM$ by $M$. After applying Fubini's theorem to (\ref{eqn:first}), we obtain

\begin{displaymath}
\int _{M} n_{T}(x,y) \, dy =  \entreceroyt \int _{S_x M} \vert \detd \vert_{\altheta} \vert \, d\theta.
\end{displaymath}

Integrating over $M$ in $x$ and using Fubini's theorem again, we get

\begin{equation}
\label{eqn:second}
\int _{M \times M} n_{T}(x,y) \, dx \, dy =  \entreceroyt \enSM \vert \detd \vert_{\altheta} \vert \, d\theta.
\end{equation}

which proves the Lemma. $\square$ \\

{\em Proof of Theorem \ref{main-teo}:} By taking $\limsup$ on both sides of (\ref{eqn:second}) and using Theorem 5.1 from Burns and Paternain \cite{burnspat}, we have proved  that

\begin{displaymath}
\htop \geq \limsup _{T \to + \infty} \frac{1}{T} \log \int _{M \times M} n_{T} (x,y) \, dx \, dy.
\end{displaymath}

If the flow admits a codimension one, continuous, invariant distribution of hyperplanes in $SM$, then the same Theorem and (\ref{eqn:second}) imply 

\begin{displaymath}
\htop = \lim _{T \to + \infty} \frac{1}{T} \log \int _{M \times M} n_{T} (x,y) \, dx \, dy. \qquad \square
\end{displaymath}

\noindent
{\bf Acknowledgement.} I am deeply grateful to Viktor Ginzburg for posing the question that led to this work and for many helpful comments and suggestions. I thank the referees for comments that helped improving the quality of exposition.

\end{document}